\documentclass[11pt]{article}

\usepackage[margin=1.1in]{geometry}
\usepackage{amssymb}

\usepackage{graphicx}
\newcommand{\bfA}{\mathbf{A}}
\newcommand{\bfb}{\mathbf{b}}
\newcommand{\x}{\mathbf{x}}

\newcommand{\epsmach}{\epsilon_{\mbox{\scriptsize mach}}}

\begin{document}

\title{Robust Numerical Tracking of One Path of a Polynomial Homotopy
      on Parallel Shared Memory Computers}
\author{Simon Telen \and Marc Van Barel\thanks{Supported by
the Research Council KU Leuven,
C1-project (Numerical Linear Algebra and Polynomial Computations),
and by
the Fund for Scientific Research--Flanders (Belgium),
G.0828.14N (Multivariate polynomial and rational interpolation and approximation),
and EOS Project no 30468160.}
\and Jan Verschelde\thanks{Supported by the National Science Foundation 
under grant DMS 1854513.}
}

\maketitle

\begin{abstract}
We consider the problem of tracking one solution path defined by
a polynomial homotopy on a parallel shared memory computer.
Our robust path tracker applies Newton's method on power series
to locate the closest singular parameter value.
On top of that, it computes singular values of the Hessians of the polynomials
in the homotopy to estimate the distance to the nearest different path.
Together, these estimates are used to compute an appropriate adaptive step size.
For $n$-dimensional problems, the cost overhead of our robust path tracker
is $O(n)$, compared to the commonly used predictor-corrector methods.
This cost overhead can be reduced
by a multithreaded program on a parallel shared memory computer.

\end{abstract}

\section{Introduction}

A polynomial homotopy is a system of polynomials in several variables
with one of the variables acting as a parameter, typically denoted by~$t$.
At $t=0$, we know the values for a solution of the system,
where the Jacobian matrix has full rank: we start at a regular solution.
With series developments we extend the values of the solution
to values of~$t > 0$.

As a demonstration of what {\em robust} in the title of this paper means,
on tracking one million paths on the 20-dimensional benchmark system
posed by Katsura~\cite{Kat90},
Table~3 of~\cite{LT09} reports 4 curve jumpings.
A curve jumping occurs when approximations from one path jump onto
another path.  In the runs with the MPI version 
for our code (reported in~\cite{TVV19})
no path failures and no curve jumpings happened.
Our path tracking algorithm applies Pad\'{e} approximants in the predictor.
These rational approximations have also been applied
to solve nonlinear systems arising in power systems~\cite{Tri12,TM16}.
In \cite{JMSW09}, Pad\'e approximants 
are used in symbolic deformation methods.

This paper describes a multithreaded version of
the robust path tracking algorithm of~\cite{TVV19}.
In~\cite{TVV19} we demonstrated the scaling of our path tracker
to polynomial homotopies with more than one million solution paths,
applying message passing for distributed memory parallel computers.
In this paper we consider shared memory parallel computers and,
starting at one single solution, we investigate the scalability
for increasing number of equations and variables,
and for an increasing number of terms in the power series developments.

As to a comparison with our MPI version used in~\cite{TVV19},
the current parallel version is made threadsafe and more efficient.
These improvements also benefit the implementation with message passing.

In addition to speedup, we ask the {\em quality up} question:
if we can afford the running time of a sequential run in double precision,
with a low degree of truncation, how many threads do we need 
(in a run which takes the same time as a sequential run)
if we want to increase the working precision and the degrees at which
we truncate the power series?

Our programming model is that of a work crew, working simultaneously 
to finish a number of jobs in a queue.
Each job in the queue is done by one single member of the work crew.
All members of the work crew have access to all data in the random
access memory of the computer.
The emphasis in this research is on the high level development 
of parallel algorithms and software~\cite{MSH11}.
The code is part of the free and open source PHCpack~\cite{Ver99},
available on github.

The parallel implementation of medium grained 
evaluation and differentiation algorithms provide good speedups.
The solution of a blocked lower triangular linear system
is most difficult to compute accurately and with good speedup.
We describe a pipelined algorithm, provide an error analysis,
and propose to apply double double and quad double arithmetic~\cite{HLB01}.

\section{Overview of the Computational Tasks}

We consider a homotopy $H$ given by $n$ polynomials $f_1, \ldots, f_n$ in $n+1$ 
variables $x_1,\ldots,x_n,t$, where $t$ is thought of as the continuation parameter. A solution path of the homotopy is denoted by $x(t)$.
For a local power series expansion $x(t) = c_0 + c_1 t + c_2 t^2 + \cdots$ of $x(t)$, where $x(t)$ is assumed analytic in a neighborhood of $t = 0$, the theorem of Fabry~\cite{Fab1896} allows us to determine
the location of the parameter value nearest to~$t=0$ where $x(t)$ is singular.
With the singular values of the Jacobian matrix $J = (\partial f_i/\partial x_j)_{1 \leq i,j \leq n}$ and the Hessian matrices of $f_1, \ldots, f_n$,
we estimate the distance to the nearest solution for $t$ fixed to zero.
The step size $\Delta t$ is the minimum of two bounds, denoted by $C$ and~$R$.

\begin{enumerate}
\item $C$ is an estimate for the nearest different solution path at $t = 0$.
      To obtain this estimate we compute the first and second partial
      derivatives at a point and organize these derivatives in the
      Jacobian and Hessian matrices.  The bound is then computed from
      the singular values of those matrices:
\begin{equation} \label{eqCbound}
   C = \frac{2\sigma_n(J)}{\sqrt{\sigma_{1,1}^2
        + \sigma_{2,1}^2 + \cdots + \sigma_{n,1}^2}},
\end{equation}
where
$\sigma_n(J)$ is the smallest singular value of the Jacobian matrix~$J$
and $\sigma_{k,1}$ is the largest singular value of the Hessian
of the $k$-th polynomial.
\item $R$ is the radius of convergence of the power series developments.
      Applying the theorem of Fabry, $R$ is computed as the ratio
      of the moduli of two consecutive coefficients in the series.
For a series truncated at degree $d$:
\begin{equation} \label{eqRbound}
   x(t) = c_0 + c_1 t + c_2 t^2 + \cdots + c_d t^d, \quad
   z = c_{d-1}/c_d, \quad R = |z|,
\end{equation}
where $z$ indicates the estimate for the location of the nearest singular parameter value.
\end{enumerate}
The computations of $R$ and $C$ require evaluation, differentiation,
and linear algebra operations.
Once $\Delta t$ is determined, the solution for the next value
of the parameter is predicted by evaluating Pad\'{e} approximants
constructed from the power series developments.
The last stage is the shift of the coefficients with $-\Delta t$,
so the next step starts again at~$t = 0$.

The stages are justified in~\cite{TVV19}.
In~\cite{TVV19}, we compared with v1.6 of Bertini~\cite{BHSW13}
(both in runs in double precision and 
in runs in adaptive precision~\cite{BHSW08})
and v1.1 of HomotopyContinuation.jl~\cite{BT18}.
In this paper we focus on parallel algorithms.

\section{Parallel Evaluation and Differentiation}

The parallel algorithms in this section are medium grained.
The jobs in the evaluation and differentiation correspond
to the polynomials in the system.
While the number of polynomials is not equal to the number of threads,
the jobs are distributed evenly among the threads.

\subsection{Algorithmic Differentiation on Power Series}

Consider a polynomial system $\bf f$ in $n$ variables
with power series (all truncated to the same fixed degree~$d$),
as coefficients;
and a vector $\bf x$ of $n$ power series, truncated to the same degree~$d$.
Our problem is to evaluate $\bf f$ at $\bf x$ and to compute all $n$
partial derivatives.
We illustrate the reverse mode of algorithmic differentiation~\cite{GW08}
with an example, on $f = x_1 x_2 x_3 x_4 x_5$.

{\small
\begin{equation} \label{eqadrevmode}
  \begin{array}{rcl}
     x_1 x_2 & = & x_1 \star x_2 \\
     x_1 x_2 x_3 & = & x_1 x_2 \star x_3 \\
     x_1 x_2 x_3 x_4 & = & x_1 x_2 x_3 \star x_4 \\
     x_1 x_2 x_3 x_4 x_5 & = & x_1 x_2 x_3 x_4 \star x_5 \\
  \end{array}
  ~
  \begin{array}{rcl}
     x_5 x_4 & = & x_5 \star x_4 \\
     x_5 x_4 x_3 & = & x_5 x_4 \star x_3 \\
     x_5 x_4 x_3 x_2 & = & x_5 x_4 x_3 \star x_2 \\
     ~ & ~ & ~ \\
  \end{array}
  ~~
  \begin{array}{rcl}
     x_1 x_3 x_4 x_5 & = & x_1 \star x_5 x_4 x_3 \\
     x_1 x_2 x_4 x_5 & = & x_1 x_2 \star x_5 x_4 \\
     x_1 x_2 x_3 x_5 & = & x_1 x_2 x_3 \star x_5 \\
     ~ & ~ & ~ \\
  \end{array}
\end{equation}
}

\noindent In the first column of~(\ref{eqadrevmode}),
we see $\frac{\partial f}{\partial x_5}$ and the evaluated $f$
on the last two rows.  The last row of the middle column gives
$\frac{\partial f}{\partial x_1}$ and the remaining partial
derivatives are in the last column of~(\ref{eqadrevmode}).

Evaluating and differentiating a product of $n$ variables in this
manner takes $3n-5$ multiplications.
For our problem, every multiplication is a convolution of two
truncated power series 
$x_i = x_{i,0} + x_{i,1} t + x_{i,2} t^2 + \cdots + x_{i,d} t^d$ and
$x_j = x_{j,0} + x_{j,1} t + x_{j,2} t^2 + \cdots + x_{j,d} t^d$,
up to degree~$d$.  Coefficients of $x_i \star x_j$ of terms higher
than~$d$ are not computed.

Any monomial is represented as the product of the variables
that occur in the monomial and the product of the monomial
divided by that product.
For example, $x_1^3 x_2 x_3^6$ is
represented as $(x_1 x_2 x_3) \cdot (x_1^2 x_3^5)$
We call the second part in this representation the common factor,
as this factor is common to all partial derivatives of the monomial.
This common factor is computed via a power table of the variables.
For every variable $x_i$,
the power table stores all powers $x_i^e$, for $e$ from 2 to the
highest occurrence in a common factor.
Once the power table is constructed, the computation of any common
factor requires at most $n-1$ multiplications
of two truncated power series.

As we expect the number of equations and variables to be a multiple
of the number of available threads, one job is the evaluation and
differentiation of one single polynomial.
Assuming each polynomial has roughly the same number of terms,
we may apply a static job scheduling mechanism.  
Let $n$ be the number of equations
(indexed from 1 to $n$), $p$ the number of threads
(labeled from 1 to $p$), where $n \geq p$.
Thread $i$ evaluates and differentiates polynomials
$i+k p$, for $k$ starting at 0, as long as $i + kp \leq n$.

\subsection{Jacobians, Hessians at a Point, and Singular Values}

If we have $n$ equations, then the computation of $C$,
defined in~(\ref{eqCbound}),
requires $n+1$ singular value decompositions, 
which can all be computed independently.

For any product of $n$ variables, after the computation of its gradient with 
the reverse mode, any element of its Hessian needs 
only a couple of multiplications, independent of~$n$.  
We illustrate this idea with an example for~$n = 8$.
The third row of the Hessian of $x_1 x_2 x_3 x_4 x_5 x_6 x_7 x_8$,
starting at the fourth column, after the zero on the diagonal is
\begin{equation}
  \begin{array}{c}
     x_1 x_2 \star x_5 x_6 x_7 x_8,
     x_1 x_2 \star x_4 \star x_6 x_7 x_8,
     x_1 x_2 x_4 \star x_5 \star x_7 x_8, \\
     x_1 x_2 x_4 x_5 \star x_6 \star x_8,
     x_1 x_2 x_4 x_5 x_6 \star x_7.
  \end{array}
\end{equation}
In the reverse mode for the gradient we already computed the forward
products $x_1 x_2$, $x_1 x_2 x_3$, $x_1 x_2 x_3 x_4$,
$x_1 x_2 x_3 x_4 x_5$, $x_1 x_2 x_3 x_4 x_5 x_6$,
and $x_1 x_2 x_3 x_4 x_5 x_6 x_7$.  
We also computed the backward products
$x_8 x_7$, $x_8 x_7 x_6$, $x_8 x_7 x_6 x_5$, $x_8 x_7 x_6 x_5 x_4$.

For a monomial $x_1^{e_1}$ $x_2^{e_2}$ $\cdots$ $x_n^{e_n}$
with higher powers $e_k > 1$, for some indices~$k$,
the off diagonal elements are multiplied with the common factor
$x_1^{e_1-1}$ $\star$ $x_2^{e_2-1}$ $\star$ $\cdots$ $\star$ $x_n^{e_n-1}$
multiplied with $e_i e_j$ at the $(i,j)$-th position in the Hessian.
The computation of this common factor requires at most $n-1$ multiplications
(fewer than $n-1$ if there are any $e_k$ equal to one),
after the computation of table which stores the values of all powers
$x_k^{e_k}$ of all values for $x_k$, for $k = 1,2,\ldots,n$.

Taking only those $m$ indices $i_k$ for which $e_{i_k} > 1$,
the common factor for all diagonal elements is
$x_{i_1}^{e_{i_1}-2}$ $x_{i_2}^{e_{i_2}-2}$ $\cdots$ 
$x_{i_m}^{e_{i_m}-2}$.
The $k$-th element on the diagonal then needs to be multiplied with
$e_{i_k}(e_{i_k}-1)$ and the product of all squares $x_{i_j}^2$,
for all $j \not= k$ for which $e_{i_j} > 1$.
The efficient computation of the sequence
$x_{i_2}^2$ $x_{i_3}^2$ $\cdots$ $x_{i_m}^2$,
$x_{i_1}^2$ $x_{i_3}^2$ $\cdots$ $x_{i_m}^2$,
$x_{i_1}^2$ $x_{i_2}^2$ $\cdots$ $x_{i_{m-1}}^2$
happens along the same lines as the computation of the gradient,
requiring $3m-5$ multiplications.

In the above paragraphs, we summarized the key ideas and results
of the application for algorithmic differentiation.
A detailed algorithmic description can be found in~\cite{Chr92}.

\section{Solving a Lower Triangular Block Linear System}

In Newton's method, the update $\Delta \x(t)$ to the power series
$\x(t)$ is computed as the solution of a linear system,
with series for the coefficient entries.

Applying linearization, we solve a sequence of as  many linear systems
(with complex numbers as coefficients), as the degree of the series.
For each linear system in the sequence, the right hand side is
computed with the solution of the previous system in the sequence.
If in each step we lose one decimal place of accuracy,
at the end of sequence we have lost as many decimal places
of accuracy as the degree of the series.

\subsection{Pipelined Solution of Matrix Series} \label{sec010}

We introduce the pipelined solution of a system of power series
by example.
Consider a power series $\bfA(t)$, with coefficients $n$-by-$n$ matrices, and
a series $\bfb(t)$, with coefficients $n$-dimensional vectors.
We want to find the solution $\x(t)$ to $\bfA(t) \x(t) = \bfb(t)$.
For series truncated to degree~5, the equation
\begin{eqnarray}
 & &  \left( A_5 t^5 + A_4 t^4 + A_3 t^3 + A_2 t^2 + A_1 t + A_0 \right)
\cdot \left( x_5 t^5 + x_4 t^4 + x_3 t^3 \right. \\
& & \left. + ~ x_2 t^2 + x_1 t + x_0 \right)
~ = ~ b_5 t^5 + b_4 t^4 + b_3 t^3 + b_2 t^2 + b_1 t + b_0
\end{eqnarray}
leads to the triangular system
(derived in~\cite{BV18} applying linearization)
\begin{eqnarray} \label{eqtriangularfirst}
   A_0 x_0 & = & b_0 \\
   A_0 x_1 & = & b_1 - A_1 x_0 \\
   A_0 x_2 & = & b_2 - A_2 x_0 - A_1 x_1 \\
   A_0 x_3 & = & b_3 - A_3 x_0 - A_2 x_1 - A_1 x_2 \\
   A_0 x_4 & = & b_4 - A_4 x_0 - A_3 x_1 - A_2 x_2 - A_1 x_3 \\
   A_0 x_5 & = & b_5 - A_5 x_0 - A_4 x_1 - A_3 x_2 - A_2 x_3 - A_1 x_4.
\label{eqtriangularlast}
\end{eqnarray}
To solve this triangular system, denote by $F_0 = F(A_0)$
the factorization of $A_0$ and $x_0 = S(F_0, b_0)$, the solution of 
$A_0 x_0 = b_0$ making use of the factorization $F_0$.
Then the equations~(\ref{eqtriangularfirst}) through~(\ref{eqtriangularlast})
are solved in the following steps.
\begin{equation} \label{eqpipelinestages}
\begin{array}{rl}
1. & F_0 = F(A_0) \\
2. & x_0 = S(F_0, b_0) \\
3. & b_1 = b_1 - A_1 x_0,~
     b_2 = b_2 - A_2 x_0,~
     b_3 = b_3 - A_3 x_0,~
     b_4 = b_4 - A_4 x_0, \\
   & b_5 = b_5 - A_5 x_0 \\
4. & x_1 = S(F_0,b_1) \\
5. & b_2 = b_2 - A_1 x_1,~
     b_3 = b_3 - A_2 x_1,~
     b_4 = b_4 - A_3 x_1,~
     b_5 = b_5 - A_4 x_1 \\
6. & x_2 = S(F_0,b_2) \\
7. & b_3 = b_3 - A_1 x_2,~
     b_4 = b_4 - A_2 x_2,~
     b_5 = b_5 - A_3 x_2 \\
8. & x_3 = S(F_0,b_3) \\
9. & b_4 = b_4 - A_1 x_3,~
     b_5 = b_5 - A_2 x_3 \\
10. & x_4 = S(F_0,b_4) \\
11. & b_5 = b_5 - A_1 x_4 \\
12. & x_5 = S(F_0,b_5) \\
\end{array}
\end{equation}
Statements on the same line can be executed simultaneously.
With 5 threads, the number of steps is reduced from 22 to 12.
For truncation degree $d$ and $d$ threads, the number of steps
in the pipelined algorithm equals $2(d+1)$.
On one thread, the number of steps equals
$2(d+1) + 1 + 2 + \cdots + d-1 = d(d-1)/2 + 2(d+1)$.
With $d$ threads, the speedup is then
\begin{equation} \label{eqpipelinedspeedup}
   \frac{d(d-1)/2 + 2(d+1)}{2(d+1)} = 1 + \frac{d(d-1)}{4 (d+1)}.
\end{equation}
As $d \rightarrow \infty$, this ratio equals $1 + d/4$.
Note that the first step is typically $O(n^3)$,
whereas the other steps are $O(n^2)$.

Observe in~(\ref{eqpipelinestages}) that the first operation
on every line is on the critical path of all possible parallel executions.
For the example in~(\ref{eqpipelinestages}) this implies that the
total number of steps will never become less than 12,
even as the number of threads goes to infinity.
The speedup of 22/12 remains the same as we reduce the number
of threads from~5 to~3, as the updates of $b_4$ and $b_5$ in step 3
can be postponed to the next step.
Likewise, the update of $b_5$ in step 5 may happen in step~6.
Generalizing this observation, the formula for the speedup
in~(\ref{eqpipelinedspeedup}) remains the same for $d/2 + 1$ threads
(instead of $d$) in case $d$ is odd.  In case $d$ is even,
then the best speedup is obtained with $d/2$ threads.

Better speedups will be obtained for finer granularities,
if the matrix factorizations are executed in parallel as well.

\subsection{Error Analysis of a Lower Triangular Block Toeplitz Solver}

In Section~\ref{sec010}, we designed a pipelined method 
to solve the following lower triangular block Toeplitz system 
of equations
\begin{equation}
\left[
\begin{array}{ccccc}
A_0 & & & & \\
A_1 & A_0 & & & \\
A_2 & A_1 & A_0 & & \\
\vdots & \vdots & \vdots & \ddots & \\
A_i & A_{i-1} & A_{i-2} & \cdots & A_0
\end{array}
\right]
\left[
\begin{array}{c}
x_0 \\ x_1 \\ x_2 \\ \vdots \\ x_i
\end{array}
\right]
=
\left[
\begin{array}{c}
b_0 \\ b_1 \\ b_2 \\ \vdots \\ b_i
\end{array}
\right].
\end{equation}
In this section, we do not intend to give a very detailed error analysis 
but indicate using a rough estimate of the norm of the blocks involved,
where and how there could be a loss of precision in some typical situations.
In our analysis we will use the Euclidean 
2-norm $\| \cdot \| = \| \cdot \|_2$ on finite dimensional complex vector 
spaces and the induced operator norm on matrices.
Without loss of generality,
we can always assume that the system is scaled such that
\begin{equation}
\| A_0 \| = \| x_0 \| = 1.
\end{equation}
Hence, assuming that the components of $x_0$ in the direction of the right 
singular vectors of $A_0$ corresponding to the larger singular values are 
not too small,
the norm of the first block $b_0$ of the right-hand side satisfies
\begin{equation}
\| b_0 \| = \| A_0 x_0 \| \lessapprox \| A_0\| \| x_0 \|. 
\end{equation}
To determine the first component $x_0$ of the solution vector,
we solve the system $A_0 x_0 = b_0$.
We solve this first system in a backward stable way, i.e.,
the computed solution $\hat{x}_0 =x_0+\Delta x_0$ can be considered
as the exact solution of the system
\begin{equation}
  A_0 \hat{x}_0  = b_0 + \Delta b_0 \qquad\mbox{ with }
  \qquad \frac{\| \Delta b_0 \|}{\| b_0 \|} \approx \epsmach.
\end{equation}
If we denote the condition number of $A_0$ by $\kappa$, we get
\begin{equation}
  \frac{\| \Delta x_0 \|}{\| x_0 \|}
  \leq \kappa \frac{\| \Delta b_0 \|}{\| b_0 \|} \leq \kappa O(\epsmach).
\end{equation}
We study now how this error influences the remainder of the calculations.
In the remaining steps, we use rough estimates of the order of magnitude 
of the different blocks $A_i$ of the coefficient matrix, the blocks $x_i$
of the solution vector and the blocks $b_i$ of the right-hand side.
First we will assume that the sizes of the blocks $x_i$ as well 
as $A_i$ behave as $\rho^i$, i.e.,
\begin{equation}
   \| x_i \| \approx \rho^i \qquad \mbox{and}\qquad \| A_i \| \approx \rho^i.
\end{equation}
Hence, also the sizes of the blocks $b_i$ behave as
\begin{equation}
   \| b_i \| \approx \rho^i .
\end{equation}
In our context, the parameter $\rho$ should be thought of as the inverse 
of the convergence radius~$R$, as defined in~(\ref{eqRbound}),
for the series expansions.  Note that when $\rho$ is larger, 
this indicates that the distance to the nearest singularity is smaller.
Consider now the second system
\begin{equation}
   A_0 x_1 = \tilde{b}_1,
\end{equation}
where $\tilde{b}_1 = b_1 - A_1x_0$.
Using the computed value $\hat{x}_0$, we find an approximation $\hat{x}_1 = x_1 + \Delta x_1$ for $x_1$ by solving the system 
\begin{equation}
  A_0 X = b_1 -  A_1 \hat{x}_0
  = b_1 - A_1x_0 - A_1 \Delta x_0 = \tilde{b}_1 - A_1 \Delta x_0.
\end{equation}
for $X$.
We have that
$\| \tilde{b}_1 \| = \|A_0 x_1 \| \approx \rho^1$.
Because $\| \Delta x_0 \| \approx \kappa \epsmach$,
this results in an absolute error $\Delta \tilde{b}_1 = - A_1 \Delta x_0$ on $\tilde{b}_1$
of size $\kappa \epsmach \rho$ or a relative error of size $\kappa \epsmach$.
Hence,
\begin{equation}
  \frac{\| \Delta x_1 \|}{\| x_1 \|} 
  \approx \kappa \frac{\| \Delta \tilde{b}_1 \|}{\| \tilde{b}_1 \|}
  \approx \kappa^2 \epsmach.
\end{equation}
In the same way, one derives that
\begin{equation}
\frac{\| \Delta x_i \|}{\|x_i\|} \approx \kappa^{i+1} \epsmach.
\end{equation}
Hence, when $\| x_i \| \approx \rho^i$ and $\|A_i\| \approx \rho^i$, we lose all precision as soon as
$\kappa^{i+1} \epsmach = O(1)$. When the matrix $A_0$ is ill-conditioned (i.e., when $\kappa$ is large), this may happen already after a few number of steps $i$.

Assuming now that $\| x_i \| \approx \rho^i$ and $\|A_i\| \approx \rho^0$,
we solve for the second block equation
\begin{equation}
   A_0 X = b_1 -  A_1 \hat{x}_0= \tilde{b}_1 - A_1 \Delta x_0
\end{equation}
with $\| \tilde{b}_1 \| = \|A_0 x_1 \| \approx \rho^1$.
However, in this case the absolute error 
$\| \Delta x_0 \| \approx \kappa \epsmach$ is not amplified
and results in an absolute error $\Delta \tilde{b}_1 = - A_1 \Delta x_0$ 
of size $\kappa \epsmach$ or a relative error of size
$\kappa \epsmach / \rho$.
If $\kappa \geq \rho$ this is the dominant error on $\tilde{b}_1$.
If $\kappa \leq \rho$, the dominant error is the error 
of computing $\tilde{b}_1$ in finite precision.
In that case, the relative error will be of size $\epsmach$.
In what follows, we'll assume that $\kappa \geq \rho$.
The other case can be treated in a similar way.
It follows that
\begin{equation}
  \frac{\| \Delta x_1 \|}{\| x_1 \|} \approx \kappa
  \frac{\| \Delta \tilde{b}_1 \|}{\| \tilde{b}_1 \|} \approx \kappa
  \frac{\kappa}{\rho} \epsmach.
\end{equation}
Next, the approximation $\hat{x}_2 = x_2 + \Delta x_2$ of $x_2$ is computed by solving
\begin{equation}
   A_0 X = b_2 - A_2 \hat{x}_0 - A_1 \hat{x}_1 = \tilde{b}_2 - A_2 \Delta x_0 - A_1 \Delta x_1
\end{equation}
for $X$, with $\tilde{b}_2 = b_2 - A_2 x_0 - A_1 x_1$ and $\| \tilde{b}_2 \| = \|A_0 x_2 \| \approx \rho^2$.
The absolute error $\Delta x_0$ plays a minor role compared to $\Delta x_1$.
The relative error on $x_1$ of magnitude $\kappa (\kappa / \rho) \epsmach$
multiplied by $A_1$ of norm $\rho$ leads to a relative error of 
magnitude $(\kappa/\rho)^2 \epsmach$ on $\tilde{b}_2$.
Hence,
\begin{equation}
  \frac{\| \Delta x_2 \|}{\| x_2 \|} \approx \kappa
  \frac{\| \Delta \tilde{b}_2 \|}{\| \tilde{b}_2 \|} \approx \kappa
  \frac{\kappa^2}{\rho^2} \epsmach.
\end{equation}
In a similar way, one derives that, when $\kappa \geq \rho$:
\begin{equation}
  \frac{\| \Delta x_i \|}{\| x_i \|} \approx  \kappa
  \frac{\kappa^i}{\rho^i} \epsmach.
\end{equation}
In an analogous way the other possibilities in the summary hereafter 
can be deduced.
Assuming that $\| x_i \| \approx \rho^i$ we have the following possibilities:
\begin{enumerate}
\item When $\| A_i \| \approx \rho^i$,
      we can not do much about the loss of accuracy:
\begin{equation}
  \frac{\| \Delta x_i \|}{\| x_i \|} \approx  \kappa^{i+1} \epsmach.
\end{equation}
\item When $\| A_i\| \approx 1^i$, we can distinguish two possibilities:
\begin{eqnarray}
  \mbox{when $\kappa \geq \rho$} &:& \frac{\| \Delta x_i \|}{\| x_i \|}
  \approx  \kappa \frac{\kappa^i}{\rho^i} \epsmach; \\
  \mbox{when $\kappa \leq \rho$} &:& \frac{\| \Delta x_i \|}{\| x_i \|}
  \approx  \kappa \epsmach.
\end{eqnarray}
The second case cannot arise when $\rho < 1$.
\end{enumerate}
We observe in computational experiments that in our path tracking method
we are usually dealing with the first case,
where $\|A_i\| \approx \rho^i$, $\|x_i \| \approx \rho^i$.
This means that the number of coefficients that we can compute 
with reasonable accuracy is bounded roughly 
by $-\log(\epsmach)/\log(\kappa)$,
where $\kappa$ is the condition number of the Jacobian $A_0$.

\subsection{Newton's Method, Rational Approximations, Coefficient Shift}

In Newton's method, the evaluation and differentiation algorithms
are followed by the solution of the matrix series system
to compute all coefficients of a power series at a regular solution
of a polynomial homotopy.  
There are two remaining stages.
Both stages use the same type of parallel algorithm,
summarized in the next two paragraphs.

A Pad\'e approximant is the quotient of two polynomials.
To construct an approximant of degree~$K$ in the numerator 
and $L$ in the denominator, 
we need the first $K+L+1$ coefficients of the power series.
Given $K$ and $L$, we truncate the power series at degree~$d = K+L$.
All components of an $n$-dimensional vector can be computed
independently from each other, so each job in the parallel algorithm
is the construction and evaluation of one Pad\'{e} approximant.

All power series are assumed to originate at $t=0$.
After incrementing the step size with $\Delta t$, we shift all
coefficients of the power series in the polynomial homotopy with $-\Delta t$,
so at the next step we start again at $t = 0$.
The shift operation happens independently for every polynomial
in the homotopy, so the threads take turns in shifting the coefficients.

As the computational experiments show, the construction
of rational approximations and the shifting of coefficients
are computationally less intensive than running Newton's method,
or than computing the Jacobian, all Hessians, and singular values
at a point. 


\section{Computational Experiments}

The goal of the computational experiments is to examine the relative
computational costs of the various stages and to detect potential
bottlenecks in the scalability.
After presenting tables for random input data,
we end with a description of a run on a cyclic $n$-root,
for $n = 64$, 96, 128,
a sample of a well known benchmark problem~\cite{Dav87}
in polynomial system solving.

Our computational experiments run
on two 22-core 2.2 GHz Intel Xeon E5-2699 processors
in a CentOS Linux workstation with 256 GB RAM.
In our speedup computation, we compare against a sequential implementation,
using the same primitive operations.

For each run on $p$ threads, we report the speedup $S(p)$,
the ratio between the serial time over the parallel execution time,
and the efficiency $E(p) = S(p)/p$.
Although our workstation has 44 cores, we stop the runs at 40 threads
to avoid measuring the interference with other unrelated processes.

The units of all times reported in the tables below are seconds
and the times themselves are elapsed wall clock times.
These times include the allocation and deallocation of all
data structures, for inputs, results, and work space.

\subsection{Random Input Data}

The randomly generated problems represent polynomial systems
of dimension~64 (or higher), with 64 (or more) terms in each 
polynomial and exponents of the variables between zero and eight.

\subsubsection{Algorithmic Differentiation on Power Series}

\begin{figure}[htb]
\begin{center}
{\includegraphics[scale=0.70]{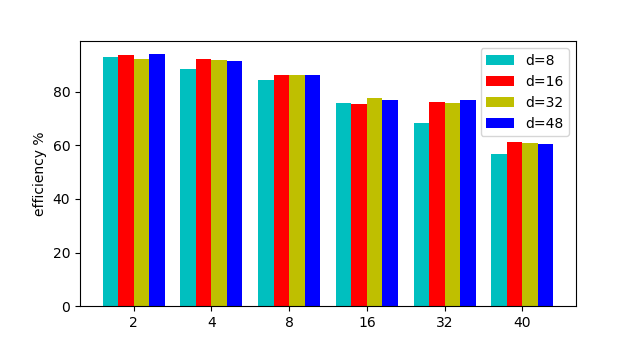}}
\caption{Efficiency plots for evaluation and differentiation
of power series, with data from Table~\ref{tabadcnv}.
Efficiency tends to decrease for increasing $p$.
The efficiency improves a little as the truncation degree $d$
of the series increases from 8, 16, 32, to~48.}
\label{figadcnv}
\end{center}
\end{figure}

The computations in Table~\ref{tabadcnv} illustrate the cost overhead
of working with power series of increasing degrees of truncation.
We start with degree $d=8$ (the default in~\cite{TVV19})
and consider the increase in wall clock times as we increase~$d$.
Reading Table~\ref{tabadcnv} diagonally, observe the quality up.
Figure~\ref{figadcnv} shows the efficiencies.

\begin{table}[hbt]
\begin{center}
\begin{tabular}{r||r|r|r||r|r|r||r|r|r||r|r|r}
    & \multicolumn{3}{c||}{$d=8$}
    & \multicolumn{3}{c||}{$d=16$}
    & \multicolumn{3}{c||}{$d=32$}
    & \multicolumn{3}{c}{$d=48$} \\
$p$~ & ~time~ & ~$S(p)$~ & $E(p)$~
     & ~time~ & ~$S(p)$~ & $E(p)$~
     & ~time~ & $S(p)$ & $E(p)$~
     & ~time~ & $S(p)$ & $E(p)$~ \\ \hline \hline
  1 & 44.851 &       &         & 154.001 &       &         & 567.731 &       &         & 1240.761 &       &         \\
  2 & 24.179 &  1.86 &  92.8\% &  82.311 &  1.87 &  93.6\% & 308.123 &  1.84 &  92.1\% &  659.332 &  1.88 &  94.1\% \\
  4 & 12.682 &  3.54 &  88.4\% &  41.782 &  3.69 &  92.2\% & 154.278 &  3.68 &  92.0\% &  339.740 &  3.65 &  91.3\% \\
  8 &  6.657 &  6.74 &  84.2\% &  22.332 &  6.90 &  86.2\% &  82.250 &  6.90 &  86.3\% &  179.424 &  6.92 &  86.4\% \\
 16 &  3.695 & 12.14 &  75.9\% &  12.747 & 12.08 &  75.5\% &  45.609 & 12.45 &  77.8\% &  100.732 & 12.32 &  76.9\% \\
 32 &  2.055 & 21.82 &  68.2\% &   6.332 & 24.32 &  76.0\% &  23.451 & 24.21 &  75.7\% &   50.428 & 24.60 &  76.9\% \\
 40 &  1.974 & 22.72 &  56.8\% &   6.303 & 24.43 &  61.1\% &  23.386 & 24.28 &  60.7\% &   51.371 & 24.15 &  60.4\%
\end{tabular}
\end{center}
\caption{Evaluation and differentiation at power series
truncated at increasing degrees~$d$, 
for increasing number of threads~$p$,
in quad double precision.}
\label{tabadcnv}
\end{table}

The drop in efficiency with $p = 40$ is because the problem size $n = 64$
is not a multiple of $p$, which results in load imbalancing.
As quad double arithmetic is already very computationally intensive,
the increase in the truncation degree $d$ does little to improve the efficiency.
Using more threads increases the memory usage, as each thread needs its own
work space for all data structures used in the computation of its gradient
with algorithmic differentiation.  In a sequential computation where gradients
are computed one after the other, there is only one vector with forward,
backward, and cross products.  When $p$ gradients are computed simultaneously,
there are $p$ work space vectors to store the intermediate forward, backward,
and cross products for each gradient.  The portion of the parallel code that
allocates and deallocates all work space vectors grows as the number of threads
increases and the wall clock times incorporate the time spent on that data
management as well.

\subsubsection{Jacobians, Hessians at a Point, and Singular Values}

Table~\ref{tabhessians} summarizes runs on the evaluation
and singular value computations on random input data,
for $n$-dimensional problems.
The $n$ polynomials have each $n$ terms, where the exponents
of the variables range from zero to eight.

\begin{table}
\begin{center}
\begin{tabular}{r|r||r|r|r||r|r|r||r|r|r}
\multicolumn{2}{c||}{~} 
  & \multicolumn{3}{c||}{double}
  & \multicolumn{3}{c||}{double double}
  & \multicolumn{3}{c}{quad double} \\
$n$~ & $p$~ & ~time~ & ~$S(p)$~ & ~$E(p)$~
     & ~time~ & ~$S(p)$~ & ~$E(p)$~
     & ~time~ & ~$S(p)$~ & ~$E(p)$~ \\ \hline \hline
64~
&   1~ & ~0.729~ &       &         & ~3.964~ &       &         & ~51.998~ &        &         \\
&   2~ & ~0.521~ & 1.40~ & ~70.0\% & ~2.329~ & 1.70~ & ~85.1\% & ~29.183~ &  1.78~ & ~89.1\% \\
&   4~ & ~0.308~ & 2.37~ & ~59.2\% & ~1.291~ & 3.07~ & ~76.8\% & ~16.458~ &  3.16~ & ~79.0\% \\
&   8~ & ~0.208~ & 3.50~ & ~43.7\% & ~0.770~ & 5.15~ & ~64.3\% &  ~9.594~ &  5.42~ & ~67.8\% \\
& ~16~ & ~0.166~ & 4.39~ & ~27.4\% & ~0.498~ & 7.96~ & ~49.8\% &  ~6.289~ &  8.27~ & ~51.7\% \\
& ~32~ & ~0.153~ & 4.77~ & ~14.9\% & ~0.406~ & 9.76~ & ~30.5\% &  ~4.692~ & 11.08~ & ~34.6\% \\
& ~40~ & ~0.129~ & 5.65~ & ~14.1\% & ~0.431~ & 9.19~ & ~23.0\% &  ~4.259~ & 12.21~ & ~30.5\% 
\\ \hline \hline
96~
&   1~ & ~3.562~ &       &         & ~18.638~ &         &          & ~240.70~ &         &         \\
&   2~ & ~2.051~ & 1.74~ & ~86.8\% & ~11.072~ &   1.68~ & ~84.17\% & ~132.76~ &   1.81~ & ~90.7\% \\
&   4~ & ~1.233~ & 2.89~ & ~72.2\% &   5.851~ &   3.19~ & ~79.64\% &   72.45~ &   3.32~ & ~83.1\% \\
&   8~ & ~0.784~ & 4.54~ & ~56.8\% &   3.374~ &   5.52~ & ~69.06\% &   41.20~ &   5.84~ & ~73.0\% \\
& ~16~ & ~0.521~ & 6.84~ & ~42.7\% &   2.188~ &   8.52~ & ~53.25\% &   25.87~ &   9.30~ & ~58.1\% \\
& ~32~ & ~0.419~ & 8.50~ & ~26.6\% &   1.612~ & ~11.56~ & ~36.13\% &   15.84~ & ~15.20~ & ~47.5\% \\
& ~40~ & ~0.398~ & 8.94~ & ~22.4\% &   1.442~ & ~12.92~ & ~32.31\% &   15.84~ & ~15.20~ & ~38.0\%
\\ \hline \hline
128~ 
&   1~ & ~12.464~ &         &         & ~62.193~ &         &         & ~730.50~ &         &         \\
&   2~ &   6.366~ &   1.96~ & ~97.9\% & ~33.213~ &   1.87~ & ~93.6\% & ~399.98~ &   1.83~ & ~91.3\% \\
&   4~ &   3.570~ &   3.49~ & ~87.3\% & ~17.436~ &   3.57~ & ~89.2\% & ~213.04~ &   3.43~ & ~85.7\% \\
&   8~ &   2.170~ &   5.75~ & ~71.8\% &   9.968~ &   6.24~ & ~78.0\% & ~119.81~ &   6.10~ & ~76.2\% \\
& ~16~ &   1.384~ &   9.01~ & ~56.3\% &   6.101~ & ~10.19~ & ~63.7\% &   73.09~ &   9.99~ & ~62.5\% \\
& ~32~ &   1.033~ & ~12.06~ & ~37.7\% &   4.138~ & ~15.03~ & ~47.9\% &   43.44~ & ~16.82~ & ~52.6\% \\
& ~40~ &   0.981~ & ~12.70~ & ~31.7\% &   3.677~ & ~16.92~ & ~42.3\% &   42.44~ & ~17.21~ & ~43.0\%
\end{tabular}
\end{center}
\caption{Evaluation of Jacobian and Hessian matrices
at a point, singular value decompositions, for $p$ threads,
in double, double double, and quad double precision.}
\label{tabhessians}
\end{table}

\begin{figure}
\begin{center}
{\includegraphics[scale=0.70]{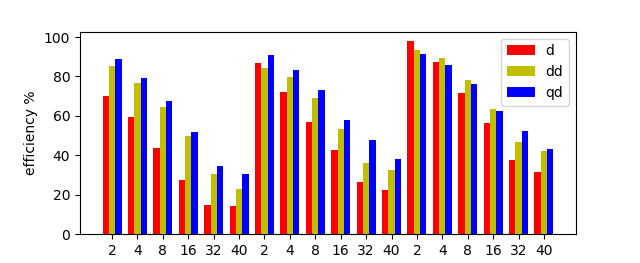}}
\caption{Efficiency plots for computing Jacobians, Hessians,
and their singular values, with data from Table~\ref{tabhessians}.
The three ranges for $p = 2, 4, 8, 16, 32, 40$ are from left to
right for $n = 64, 96$, and 128 respectively.
Efficiency decreases for increasing values of $p$.
Efficiency increases for increasing values of $n$
and for increased precision, 
where d = double, dd = double double, and qd = quad double.}
\label{fighessians}
\end{center}
\end{figure}

Reading the columns of Table~\ref{tabhessians} vertically,
we observe increasing speedups, which increase as~$n$ increases.
Reading Table~\ref{tabhessians} horizontally, 
we observe the cost overhead of the arithmetic.
To see how many threads are needed to compensate for this overhead,
read Table~\ref{tabhessians} diagonally.
Figure~\ref{fighessians} shows the efficiencies.

To explain the drop in efficiencies we apply the same reasoning
as before and point out that the work space increases even more
as more threads are applied, because the total memory consumption
has increased with the two dimensional Hessian matrices.

\subsubsection{Pipelined Solution of Matrix Series}

\begin{figure}[b!]
\begin{center}
{\includegraphics[scale=0.70]{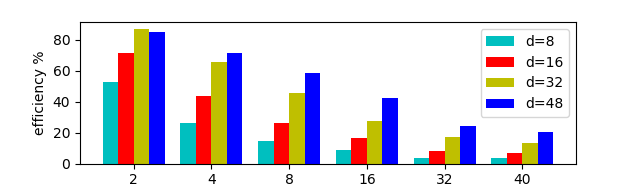}}
\caption{Efficiency plots for pipelined solution of a matrix series
with data from Table~\ref{tabserlin}.
Efficiency tends to decrease for increasing $p$
and increase for increasing $d$.}
\label{figpipe}
\end{center}
\end{figure}

Elapsed wall clock times and speedups are listed
in Table~\ref{tabserlin}, on randomly generated linear systems
of 64 equations in 64 unknowns, for series truncated to increasing degrees.
The dimensions are consistent with the setup of Table~\ref{tabadcnv},
to relate the cost of linear system solving to the cost of evaluation
and differentiations.
Figure~\ref{figpipe} shows the efficiencies.

\begin{table}[hbt]
\begin{center}
\begin{tabular}{r||r|r|r||r|r|r||r|r|r||r|r|r}
    & \multicolumn{3}{c||}{$d=8$}
    & \multicolumn{3}{c||}{$d=16$}
    & \multicolumn{3}{c||}{$d=32$}
    & \multicolumn{3}{c}{$d=48$} \\
$p$~ & ~time~ & $S(p)$ & $E(p)$~
     & ~time~ & $S(p)$ & $E(p)$~
     & ~time~ & $S(p)$ & $E(p)$~
     & ~time~ & $S(p)$ & $E(p)$~ \\ \hline \hline
  1 & 0.232 &      &         & 0.605 &      &         & 2.022 &      &         & 4.322        &         \\
  2 & 0.222 & 1.05 &  52.4\% & 0.422 & 1.44 &  71.7\% & 1.162 & 1.74 &  87.0\% & 2.553 & 1.69 &  84.7\% \\
  4 & 0.218 & 1.07 &  26.6\% & 0.349 & 1.74 &  43.4\% & 0.775 & 2.61 &  65.3\% & 1.512 & 2.86 &  71.5\% \\
  8 & 0.198 & 1.18 &  14.7\% & 0.291 & 2.08 &  26.0\% & 0.554 & 3.65 &  45.6\% & 0.927 & 4.66 &  58.3\% \\
 16 & 0.166 & 1.40 &   8.7\% & 0.225 & 2.69 &  16.8\% & 0.461 & 4.39 &  27.5\% & 0.636 & 6.80 &  42.5\% \\
 32 & 0.197 & 1.18 &   3.7\% & 0.225 & 2.69 &   8.4\% & 0.371 & 5.45 &  17.0\% & 0.554 & 7.81 &  24.4\% \\
 40 & 0.166 & 1.40 &   3.5\% & 0.227 & 2.67 &   6.7\% & 0.369 & 5.48 &  13.7\% & 0.531 & 8.14 &  20.3\%
\end{tabular}
\end{center}
\caption{Solving a linear system for power series
truncated at increasing degrees~$d$, 
for increasing number of threads~$p$,
in quad double precision.}
\label{tabserlin}
\end{table}

Consistent with the above analysis,
the speedups in Table~\ref{tabserlin} level off for $p > d/2$.
A diagonal reading shows that with multithreading,
we can keep the time below one second, while increasing the degree
of the truncation from 8 to~48.
Relative to the cost of evaluation and differentiation,
the seconds in Table~\ref{tabserlin} are significantly smaller than
the seconds in Table~\ref{tabadcnv}.

\subsubsection{Multithreaded Newton's Method on Power Series}

\begin{figure}[hbt]
\begin{center}
{\includegraphics[scale=0.70]{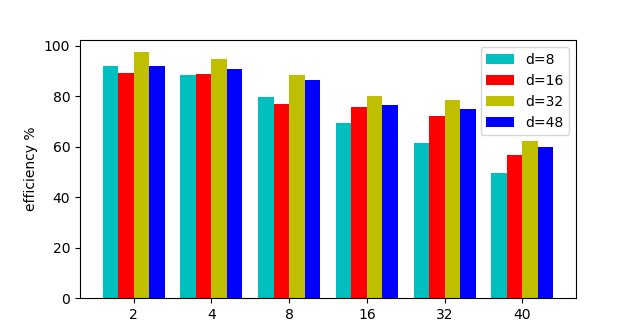}}
\caption{Efficiency plots for running Newton's method
with data from Table~\ref{tabnewton}.
Efficiency tends to decrease for increasing $p$
and increase for increasing degree $d$.}
\label{fignewton}
\end{center}
\end{figure}

In the randomly generated problems, we add the parameter $t$
to every polynomial to obtain a Newton homotopy.
The elapsed wall clock times in Table~\ref{tabnewton} 
come from running Newton's method, which requires the repeated
evaluation, differentiation, and linear system solving.
The dimensions of the randomly generated problems are
64 equations in 64 variables, with 8 as the highest degree 
in each variable.  The parameter~$t$ appears with degree one.
Figure~\ref{fignewton} shows the efficiencies.

\begin{table}[htb]
\begin{center}
\begin{tabular}{r||r|r|r||r|r|r||r|r|r||r|r|r}
    & \multicolumn{3}{c||}{$d=8$}
    & \multicolumn{3}{c||}{$d=16$}
    & \multicolumn{3}{c||}{$d=32$}
    & \multicolumn{3}{c}{$d=48$} \\
$p$~ & ~time~ & $S(p)$ & $E(p)$~
     & ~time~ & $S(p)$ & $E(p)$~
     & ~time~ & $S(p)$ & $E(p)$~
     & ~time~ & $S(p)$ & $E(p)$~ \\ \hline \hline
  1 & 347.854 &       &         & 1176.887 &       &         & 4525.080 &       &         & 7005.914 &       &         \\
  2 & 188.922 &  1.84 &  92.1\% &  658.935 &  1.79 &  89.3\% & 2323.203 &  1.95 &  97.4\% & 3806.198 &  1.84 &  92.0\% \\
  4 &  98.281 &  3.54 &  88.5\% &  330.497 &  3.56 &  89.0\% & 1193.762 &  3.79 &  94.8\% & 1925.040 &  3.64 &  91.0\% \\
  8 &  54.551 &  6.38 &  79.7\% &  191.575 &  6.14 &  76.8\% &  638.208 &  7.09 &  88.6\% & 1014.856 &  6.90 &  86.3\% \\
 16 &  31.262 & 11.13 &  69.5\% &   97.342 & 12.09 &  75.6\% &  352.103 & 12.85 &  80.3\% &  571.258 & 12.26 &  76.7\% \\
 32 &  17.624 & 19.74 &  61.7\% &   50.809 & 23.16 &  72.4\% &  180.318 & 25.60 &  78.4\% &  291.923 & 24.00 &  75.0\% \\
 40 &  17.456 & 19.93 &  49.8\% &   51.701 & 22.76 &  56.9\% &  181.563 & 24.92 &  62.3\% &  292.552 & 23.95 &  59.9\%
\end{tabular}
\end{center}
\caption{Running 8 steps with Newton's method for power series
truncated at increasing degrees~$d$, 
for increasing number of threads~$p$,
in quad double precision.}
\label{tabnewton}
\end{table}

The improvement in the efficiencies as the degrees increase
can be explained by the improvement in the efficiencies
in the pipelined solution of matrix series, see Figure~\ref{figpipe}.

\subsubsection{Rational Approximations}

\begin{figure}[b!]
\begin{center}
{\includegraphics[scale=0.70]{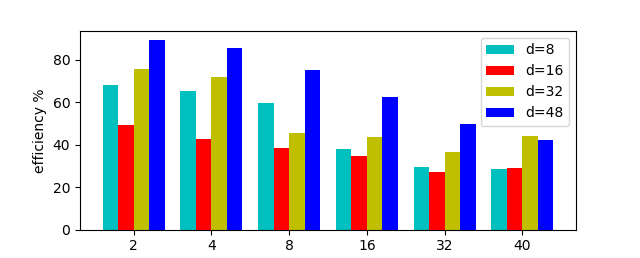}}
\caption{Efficiency plots for rational approximations
with data from Table~\ref{tabpade}.}
\label{figpade}
\end{center}
\end{figure}

In Table~\ref{tabpade}, wall clock times and speedups are listed
for the construction and evaluation of vectors of Pad\'{e} approximants,
of dimension~64 and for increasing degrees $d = 8, 16, 24$, and~$32$.
For each $d$, we take $K = L = d/2$.
Figure~\ref{figpade} shows the efficiencies.
The fast drop in efficiency for $d=8$ is due to the tiny wall clock times.
There is not much that can be improved with multithreading once the time
drops below 10 milliseconds.

\begin{table}[t!]
\begin{center}
\begin{tabular}{r||r|r|r||r|r|r||r|r|r||r|r|r}
    & \multicolumn{3}{c||}{$d=8$}
    & \multicolumn{3}{c||}{$d=16$}
    & \multicolumn{3}{c||}{$d=32$}
    & \multicolumn{3}{c}{$d=48$} \\
$p$~ & ~time~ & ~$S(p)$~ & ~$E(p)$~
     & ~time~ & ~$S(p)$~ & ~$E(p)$~
     & ~time~ & ~$S(p)$~ & ~$E(p)$~
     & ~time~ & ~$S(p)$~ & ~$E(p)$~ \\ \hline \hline
  1 & 0.034 &       &         & 0.109 &       &         & 0.684 &       &         & 2.193 &       &         \\
  2 & 0.025 &  1.36 &  68.1\% & 0.110 &  0.99 &  49.4\% & 0.452 &  1.51 &  75.6\% & 1.231 &  1.78 &  89.1\% \\
  4 & 0.013 &  2.61 &  65.2\% & 0.064 &  1.71 &  42.6\% & 0.238 &  2.87 &  71.8\% & 0.642 &  3.42 &  85.4\% \\
  8 & 0.007 &  4.79 &  59.8\% & 0.035 &  3.07 &  38.4\% & 0.189 &  3.63 &  45.4\% & 0.365 &  6.01 &  75.1\% \\
 16 & 0.006 &  6.09 &  38.1\% & 0.020 &  5.52 &  34.5\% & 0.098 &  6.96 &  43.5\% & 0.219 & 10.00 &  62.5\% \\
 32 & 0.004 &  9.47 &  29.6\% & 0.013 &  8.66 &  27.1\% & 0.058 & 11.70 &  36.6\% & 0.138 & 15.89 &  49.7\% \\
 40 & 0.003 & 11.48 &  28.7\% & 0.009 & 11.57 &  28.9\% & 0.039 & 17.58 &  43.9\% & 0.130 & 16.93 &  42.3\%
\end{tabular}
\end{center}
\caption{Construction and evaluation of Pad\'{e} approximants
for increasing degrees~$d$, 
for increasing number of threads~$p$,
in quad double precision.}
\label{tabpade}
\end{table}

\subsubsection{Shifting the Coefficients of the Power Series}

\begin{figure}[b!]
\begin{center}
{\includegraphics[scale=0.70]{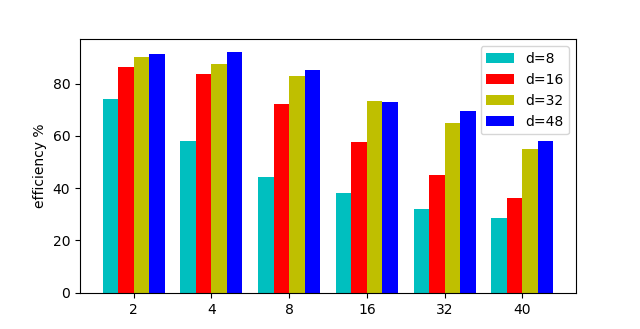}}
\caption{Efficiency plots for shifting series of a polynomial homotopy
with data from Table~\ref{tabshift}.
Efficiency tends to decrease for increasing $p$
and increase for increasing $d$.}
\label{figshift}
\end{center}
\end{figure}

Table~\ref{tabshift} summarizes experiments on a randomly generated system
of 64 polynomials in 64 unknowns, with 64 terms in every polynomial.
Figure~\ref{figshift} shows the efficiencies.

\begin{table}[t!]
\begin{center}
\begin{tabular}{r||r|r|r||r|r|r||r|r|r||r|r|r}
    & \multicolumn{3}{c||}{$d=8$}
    & \multicolumn{3}{c||}{$d=16$}
    & \multicolumn{3}{c||}{$d=32$}
    & \multicolumn{3}{c}{$d=48$} \\
$p$~ & ~time~ & $S(p)$ & $E(p)$~
     & ~time~ & $S(p)$ & $E(p)$~
     & ~time~ & $S(p)$ & $E(p)$~
     & ~time~ & $S(p)$ & $E(p)$~ \\ \hline \hline
  1 & 0.358 &       &         & 1.667 &       &         & 9.248 &       &         & 26.906 &       &         \\
  2 & 0.242 &  1.48 &  74.0\% & 0.964 &  1.73 &  86.5\% & 5.134 &  1.80 &  90.1\% & 14.718 &  1.83 &  91.4\% \\
  4 & 0.154 &  2.32 &  58.0\% & 0.498 &  3.35 &  83.8\% & 2.642 &  3.50 &  87.5\% &  7.294 &  3.69 &  92.2\% \\
  8 & 0.101 &  3.55 &  44.4\% & 0.289 &  5.77 &  72.1\% & 1.392 &  6.64 &  83.0\% &  3.941 &  6.83 &  85.3\% \\
 16 & 0.058 &  6.13 &  38.3\% & 0.181 &  9.23 &  57.7\% & 0.788 & 11.73 &  73.3\% &  2.307 & 11.66 &  72.9\% \\
 32 & 0.035 & 10.30 &  32.2\% & 0.116 & 14.40 &  45.0\% & 0.445 & 20.80 &  65.0\% &  1.212 & 22.20 &  69.4\% \\
 40 & 0.031 & 11.49 &  28.7\% & 0.115 & 14.51 &  36.3\% & 0.419 & 22.05 &  55.1\% &  1.156 & 23.28 &  58.2\%
\end{tabular}
\end{center}
\caption{Shifting the coefficients of a polynomial homotopy,
for increasing degrees~$d$, 
for increasing number of threads~$p$,
in quad double precision.}
\label{tabshift}
\end{table}

\subsubsection{Proportional Costs}

Comparing the times in Tables~\ref{tabadcnv}, \ref{tabhessians},
\ref{tabserlin}, \ref{tabpade}, and~\ref{tabshift},
we get an impression on the relative costs of the different tasks.
The evaluation and differentiation at power series,
truncated at $d = 8$ dominates the cost with 348 seconds for one thread,
or 17 seconds for 40 threads, in quad double arithmetic, 
from Table~\ref{tabadcnv}.  The second largest cost comes from
Table~\ref{tabhessians}, for $n=64$, in quad double arithmetic:
52 seconds for one thread, or 4 seconds on 40 threads.
The other three stages take less than one second on one thread.

\subsection{One Cyclic $n$-Root, $n = 64, 96, 128$}

Our algorithms are developed to run on highly nonlinear problems
such as the cyclic $n$-roots problem:
\begin{equation} \label{eqcyclicsys}
  \left\{
   \begin{array}{c}
   x_{0}+x_{1}+ \cdots +x_{n-1}=0 \\
   i = 2, 4, \ldots, n-1: 
    \displaystyle\sum_{j=0}^{n-1} ~ \prod_{k=j}^{j+i-1}
    x_{k~{\rm mod}~n}=0 \\
   x_{0}x_{1}x_{2} \cdots x_{n-1} - 1 = 0. \\
  \end{array}
  \right.
\end{equation}
This well known benchmark problem in polynomial system solving
is important in the study of biunimodular vectors~\cite{FR15}.

\subsubsection{Problem Setup}

By Backelin's Lemma~\cite{Bac89}, we know there is a 7-dimensional
surface of cyclic 64-roots, along with a recipe to generate points
on this surface.  
To generate points, a tropical formulation of 
Backelin's Lemma~\cite{AV13} is used.
The surface has degree eight.
Seven linear equations with random complex coefficients
are added to obtain isolated points on the surface.
The addition of seven linear equations gives 71 equations in 64 variables.
As in~\cite{SV00}, we add extra slack variables in an embedding to
obtain an equivalent square 71-dimensional system.
Similary, there is a 3-dimensional surface of cyclic 96-roots
and again a 7-dimensional surface of cyclic 128-roots.

In~\cite{VY15a}, running the typical predictor-corrector methods,
we experienced that the hardware double precision is no longer
sufficient to track a solution path on this 7-dimensional surface
of cyclic 64-roots.  Observe the high degrees of the polynomials
in~(\ref{eqcyclicsys}).

\begin{figure}[hbt]
\begin{center}
{\includegraphics[scale=0.70]{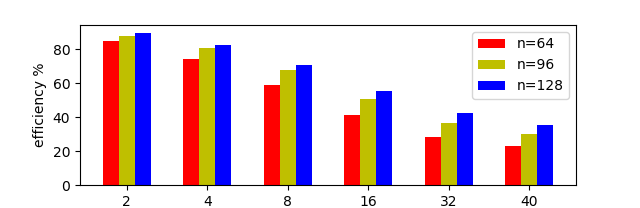}}
\caption{Efficiency plots for computing $C$ for one cyclic $n$-root,
for $n = 64$, 96, 128, with data from Table~\ref{tabcycCbound}.
Efficiency decreases for increasing~$p$
and increases for increasing~$n$.}
\label{figcurvaturecyclic}
\end{center}
\end{figure}
\begin{figure}[hbt]
\begin{center}
{\includegraphics[height=6cm]{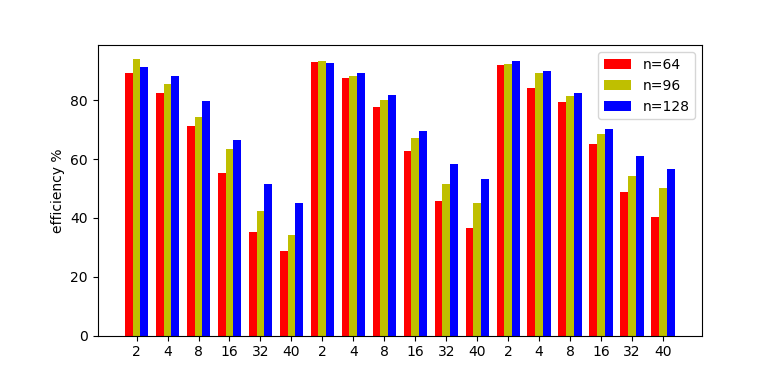}}
\caption{Efficiency plots for computing $R$ for one cyclic $n$-root,
for $n = 64$, 96, 128, for degrees $d = 8$, 16, 24, 
with data from Table~\ref{tabcycCbound}.
Efficiency decreases for increasing $p$.
Efficiency increases as $n$ and/or $d$ increase.}
\label{figradiuscyclic}
\end{center}
\end{figure}

Table~\ref{tabcycCbound} contains wall clock times, speedups
and efficiencies for computing the curvature bound $C$ for one
cyclic $n$-root.  Efficiencies are shown in Figure~\ref{figcurvaturecyclic}.
Table~\ref{tabcycRbound} contains wall clock times, speedups
and efficiencies for computing the radius bound $R$ for one
cyclic $n$-root.  See Figure~\ref{figradiuscyclic}.

\begin{table}[hbt]
\begin{center}
\begin{tabular}{r||r|r|r||r|r|r||r|r|r}
    & \multicolumn{3}{c||}{$n=64$}
    & \multicolumn{3}{c||}{$n=96$}
    & \multicolumn{3}{c}{$n=128$} \\
$p$ & \multicolumn{1}{c|}{time} & $S(p)$ & $E(p)$~
    & \multicolumn{1}{c|}{time} & $S(p)$ & $E(p)$~
    & \multicolumn{1}{c|}{time} & $S(p)$ & $E(p)$~ \\ \hline \hline
  1 & 36.862 &      &         & 152.457 &       &         & 471.719 &       &         \\
  2 & 21.765 & 1.69 &  84.7\% &  87.171 &  1.75 &  87.5\% & 262.678 &  1.80 &  89.8\% \\
  4 & 12.390 & 2.98 &  74.4\% &  47.268 &  3.23 &  80.6\% & 143.262 &  3.29 &  82.3\% \\
  8 &  7.797 & 4.73 &  59.1\% &  28.127 &  5.42 &  67.8\% &  83.044 &  5.68 &  71.0\% \\
 16 &  5.600 & 6.58 &  41.1\% &  18.772 &  8.12 &  50.8\% &  53.235 &  8.86 &  55.4\% \\
 32 &  4.059 & 9.08 &  28.4\% &  12.988 & 11.74 &  36.7\% &  34.800 & 13.56 &  42.4\% \\
 40 &  4.046 & 9.11 &  22.8\% &  12.760 & 11.95 &  29.9\% &  33.645 & 14.02 &  35.1\%
\end{tabular}
\end{center}
\caption{Computing $C$ for one cyclic $n$-root,
for $n = 64$, 96, 128,
for an increasing number of threads~$p$,
in quad double precision.}
\label{tabcycCbound}
\end{table}

\begin{table}[h!]
\begin{center}
\begin{tabular}{r|r||r|r|r||r|r|r||r|r|r}
\multicolumn{2}{c|}{~} & \multicolumn{3}{c||}{$n=64$}
                       & \multicolumn{3}{c||}{$n=96$}
                       & \multicolumn{3}{c}{$n=128$} \\
$d$ & $p$ & \multicolumn{1}{c|}{time} & $S(p)$ & $E(p)$~
          & \multicolumn{1}{c|}{time} & $S(p)$ & $E(p)$~
          & \multicolumn{1}{c|}{time} & $S(p)$ & $E(p)$~ \\ \hline \hline
8 &  1 & 139.185 &       &         & 483.137 &       &         & 1123.020 &       &         \\
  &  2 &  78.057 &  1.78 &  89.2\% & 257.023 &  1.88 &  94.0\% &  614.750 &  1.83 &  91.3\% \\
  &  4 &  42.106 &  3.31 &  82.6\% & 141.329 &  3.42 &  85.5\% &  318.129 &  3.53 &  88.3\% \\
  &  8 &  24.452 &  5.69 &  71.2\% &  81.308 &  5.94 &  74.3\% &  176.408 &  6.37 &  79.6\% \\
  & 16 &  15.716 &  8.86 &  55.4\% &  47.585 & 10.15 &  63.5\% &  105.747 & 10.62 &  66.4\% \\
  & 32 &  12.370 & 11.25 &  35.2\% &  35.529 & 13.60 &  42.5\% &   68.025 & 16.51 &  51.6\% \\
  & 40 &  12.084 & 11.52 &  28.8\% &  35.212 & 13.72 &  34.3\% &   62.119 & 18.08 &  45.2\%
\\ \hline
16 &  1 & 477.956 &       &         & 1606.174 &       &         & 3829.567 &       &         \\
   &  2 & 256.846 &  1.86 &  93.0\% &  861.214 &  1.87 &  93.3\% & 2066.680 &  1.85 &  92.7\% \\
   &  4 & 136.731 &  3.50 &  87.4\% &  454.917 &  3.53 &  88.3\% & 1072.106 &  3.57 &  89.3\% \\
   &  8 &  77.034 &  6.20 &  77.6\% &  251.066 &  6.40 &  80.0\% &  584.905 &  6.55 &  81.8\% \\
   & 16 &  47.473 & 10.07 &  62.9\% &  149.288 & 10.76 &  67.2\% &  344.430 & 11.12 &  69.5\% \\
   & 32 &  32.744 & 14.60 &  45.6\% &   97.514 & 16.47 &  51.5\% &  205.034 & 18.68 &  58.4\% \\
   & 40 &  32.869 & 14.54 &  36.4\% &   89.260 & 18.00 &  45.0\% &  180.207 & 21.25 &  53.1\%
\\ \hline
24 &  1 & 1023.968 &       &         & 3420.576 &       &         & 8146.102 &       &         \\
   &  2 &  555.771 &  1.84 &  92.1\% & 1855.748 &  1.84 &  92.2\% & 4360.870 &  1.87 &  93.4\% \\
   &  4 &  304.480 &  3.36 &  84.1\% &  956.443 &  3.58 &  89.4\% & 2268.632 &  3.59 &  89.8\% \\
   &  8 &  160.978 &  6.36 &  79.5\% &  523.763 &  6.53 &  81.6\% & 1235.338 &  6.59 &  82.4\% \\
   & 16 &   98.336 & 10.41 &  65.1\% &  312.698 & 10.94 &  68.4\% &  726.287 & 11.22 &  70.1\% \\
   & 32 &   65.448 & 15.65 &  48.9\% &  196.488 & 17.41 &  54.4\% &  416.735 & 19.55 &  61.1\% \\
   & 40 &   63.412 & 16.15 &  40.4\% &  170.474 & 20.07 &  50.2\% &  360.419 & 22.60 &  56.5\%
\end{tabular}
\end{center}
\caption{Computing $R$ for one cyclic $n$-root,
for $n = 64$, 96, 128,
for degrees~$d = 8, 16, 24$, and
for an increasing number of threads~$p$,
in quad double precision.}
\label{tabcycRbound}
\end{table}

For $n=64$, the inverse condition number of the Jacobian matrix is
estimated as $3.9\mathrm{E}{-5}$
and after 8 iterations, the maximum norm of the last vector in the 
last update to the series equals respectively $4.6\mathrm{E}{-44}$,
$1.1\mathrm{E}{-24}$, and $4.1\mathrm{E}{-5}$, for $d = 8, 16$, and~24.
For $n=96$, the estimated inverse condition number is $2.0\mathrm{E}{-4}$ 
and the maximum norm for $d = 8, 16$, and~24 is then respectively $
1.4\mathrm{E}{-47}$, $9.6\mathrm{E}{-31}$,
and~$7.3\mathrm{E}{-14}$.  The condition worsens for $n=128$,
estimated at $4.6\mathrm{E}{-6}$ and then for $d=8$, the maximum norm
of the last update vector is $2.2\mathrm{E}{-30}$.
For $d=16$ and~24, the largest maximum norm less than one occurs 
at the coefficients with $t^{15}$ and equals about $1.1\mathrm{E}{-1}$.

\section{Conclusions}

The cost overhead of our robust path tracker is $O(n)$, compared with 
the current numerical predictor-corrector algorithms.
For $n = 64$, we expect a cost overhead factor of about~64.  
We interpret the speedups in Table~\ref{tabcycCbound}
and Table~\ref{tabcycRbound} as follows.
With a speedup of about 10, then this factor drops to about~6. 

The plan is to integrate the new algorithms in the
parallel blackbox solver~\cite{Ver18}.

\bibliographystyle{splncs04}

\begin{thebibliography}{19}

\bibitem{AV13}
D.~Adrovic and J.~Verschelde.
\newblock Polyhedral methods for space curves exploiting symmetry applied to
  the cyclic $n$-roots problem.
\newblock 
  In {\em Computer Algebra in Scientific Computing, 15th International
  Workshop, CASC 2013}, pages 10--29. Springer-Verlag, 2013.

\bibitem{Bac89}
J.~Backelin.
\newblock Square multiples n give infinitely many cyclic n-roots.
\newblock Reports, Matematiska Institutionen~8, Stockholms universitet, 1989.

\bibitem{BHSW08}
D.~J. Bates, J.~D. Hauenstein, A.~J. Sommese, and C.~W. Wampler.
\newblock Adaptive multiprecision path tracking.
\newblock {\em SIAM J.\ Numer.\ Anal.}, 46(2):722--746, 2008.

\bibitem{BHSW13}
D.~J. Bates, J.~D. Hauenstein, A.~J. Sommese, and C.~W. Wampler.
\newblock {\em Numerically Solving Polynomial Systems with {B}ertini},
  volume~25.
\newblock SIAM, 2013.

\bibitem{BV18}
N.~Bliss and J.~Verschelde.
\newblock The method of {G}auss--{N}ewton to compute power series solutions of
  polynomial homotopies.
\newblock {\em Linear Algebra Appl.}, 542:569--588, 2018.

\bibitem{BT18}
P.~Breiding and S.~Timme.
\newblock Homotopy{C}ontinuation.jl: A package for homotopy continuation in
  {J}ulia.
\newblock In {\em International Congress on Mathematical Software}, pages
  458--465. Springer, 2018.

\bibitem{Chr92}
B. Christianson.
\newblock Automatic {H}essians by reverse accumulation.
\newblock {\em IMA Journal of Numerical Analysis} 12:135--150, 1992.

\bibitem{Dav87}
J. H. Davenport.
\newblock Looking at a set of equations.
\newblock Bath Computer Science Technical Report 87-06, 1987.

\bibitem{Fab1896}
E.~Fabry.
\newblock Sur les points singuliers d'une fonction donn{\'e}e par son
  d{\'e}veloppement en s{\'e}rie et l'impossibilit{\'e} du prolongement
  analytique dans des cas tr{\`e}s g{\'e}n{\'e}raux.
\newblock In {\em Annales scientifiques de l'{\'E}cole Normale Sup{\'e}rieure},
  volume~13, pages 367--399. Elsevier, 1896.

\bibitem{FR15}
H.~F{\"{u}}hr and Z.~Rzeszotnik.
\newblock On biunimodular vectors for unitary matrices.
\newblock {\em Linear Algebra Appl.}, 484:86--129, 2015.

\bibitem{GW08}
A.~Griewank and A.~Walther.
\newblock {\em Evaluating derivatives: principles and techniques of algorithmic
  differentiation}, volume 105.
\newblock SIAM, 2008.

\bibitem{HLB01}
Y. Hida, X. S. Li, and D. H. Bailey.
\newblock Algorithms for quad-double precision floating point arithmetic.
\newblock In the {\em Proceedings  of the 15th IEEE Symposium on Computer
Arithmetic (Arith-15 2001)}, pages 155--162. IEEE Computer Society,
2001. 

\bibitem{JMSW09}
G.~Jeronimo, G.~Matera, P.~Solern{\'{o}}, and A.~Waissbein.
\newblock Deformation techniques for sparse systems.
\newblock {\em Found.\ Comput.\ Math.}, 9:1--50, 2009.

\bibitem{Kat90}
S.~Katsura.
\newblock Spin glass problem by the method of integral equation of the
  effective field.
\newblock In M.~Coutinho-Filho and S.~Resende, editors, {\em New Trends in
  Magnetism}, pages 110--121. World Scientific, 1990.

\bibitem{LT09}
T.~Li and C.~Tsai.
\newblock {HOM4PS-2.0para}: {P}arallelization of {HOM4PS-2.0} for solving
  polynomial systems.
\newblock {\em Parallel Computing}, 35(4):226--238, 2009.

\bibitem{MSH11}
J. W. McCormick, F. Singhoff, and J. Hugues.
\newblock {\em Building Parallel, Embedded, and Real-Time Applications 
with Ada.} 
\newblock Cambridge University Press, 2011.

\bibitem{SV00}
A.J. Sommese and J.~Verschelde.
\newblock Numerical homotopies to compute generic points on positive
  dimensional algebraic sets.
\newblock {\em J.\ of Complexity}, 16(3):572--602, 2000.

\bibitem{TVV19}
S. Telen, M. Van Barel, and J. Verschelde.
\newblock A Robust Numerical Path Tracking Algorithm 
          for Polynomial Homotopy Continuation.
\newblock {\tt arXiv:1909.04984}.

\bibitem{Tri12}
A.~Trias.
\newblock The holomorphic embedding load flow method.
\newblock In {\em 2012 IEEE Power and Energy Society General Meeting}, pages
  1--8. IEEE, 2012.

\bibitem{TM16}
A.~Trias and J.~L. Martin.
\newblock The holomorphic embedding loadflow method for {DC} power systems and
  nonlinear {DC} circuits.
\newblock {\em IEEE Transactions on Circuits and Systems}, 63(2):322--333,
  2016.

\bibitem{Ver99}
J.~Verschelde.
\newblock Algorithm 795: {PHCpack}: A general-purpose solver for polynomial
  systems by homotopy continuation.
\newblock {\em ACM Transactions on Mathematical Software (TOMS)},
  25(2):251--276, 1999.

\bibitem{Ver18}
J.~Verschelde.
\newblock A blackbox polynomial system solver on parallel shared memory
  computers.
\newblock 
  In {\em Computer Algebra in Scientific Computing, 20th International
  Workshop, CASC 2018}, pages 361--375. Springer-Verlag, 2018.

\bibitem{VY15a}
J.~Verschelde and X.~Yu.
\newblock Accelerating polynomial homotopy continuation on a graphics
  processing unit with double double and quad double arithmetic.
\newblock In {\em Proceedings of the
  7th International Workshop on Parallel Symbolic Computation (PASCO 2015)},
  pages 109--118. ACM, 2015.

\end{thebibliography}

\end{document}